\documentclass{amsart}
\usepackage[latin9]{inputenc}
\usepackage{mathtools}
\usepackage{amsthm}
\usepackage{amssymb}
\usepackage{esint}
\usepackage{enumitem}
 
\makeatletter

 \newtheorem{thm}{Theorem}[section]

\renewcommand{\P}{\mu}
\renewcommand{\phi}{\varphi}

\begin{document}

\title{Some recent developments in quantization of fractal measures}

\author{Marc Kesseb\"ohmer}

\address{Fachbereich 3 -- Mathematik und Informatik\\
  Universit\"at Bremen\\ Bibliothekstr. 1\\Bremen 28359\\Germany}

\email{mhk@math.uni-bremen.de}

\author{Sanguo Zhu}

\address{School of Mathematics and Physics\\
Jiangsu University
of Technology
\\ Changzhou 213001\\China}

\email{sgzhu@jsut.edu.cn}
\thanks{S. Zhu was supported by China Scholarship Council No. 201308320049}
\date{January 6, 2015}
\begin{abstract}
We give an overview on the quantization problem for fractal measures,
including some related results and methods which have been developed
in the last decades. Based on the work of Graf and Luschgy, we propose
a three-step procedure to estimate the quantization errors. We survey
some recent progress, which makes use of this procedure, including
the quantization for self-affine measures, Markov-type measures on
graph-directed fractals, and product measures on multiscale Moran
sets. Several open problems are mentioned.  
\end{abstract}



\maketitle

\section{Introduction}

The quantization problem for probability measures originated in information
theory and certain areas of engineering technology such as image compression and
data processing. In the past decades, this problem has been rigorously
studied by mathematicians and the field of quantization theory emerged.
Recently, this theory has also been found to have promising applications
in numerical integrations and mathematical finance (see e.g. \cite{SpaceQuantizationMethodIntegration,PagesWilbergt,SwingAmericanDualPrimal}).
Mathematically we are concerned with the asymptotics of the errors
in the approximation of a given probability measure with finitely
supported probability measures in the sense of $L_{r}$-metrics. More
precisely, for every $n\in\mathbb{N}$, we set $\mathcal{D}_{n}\coloneqq\{\alpha\subset\mathbb{R}^{q}:1\leq{\rm card}(\alpha)\leq n\}$.
Let $\P$ be a Borel probability measure on $\mathbb{R}^{q}$, $q\in\mathbb{N}$,
and let $r\in[0,\infty)$. The $n$-th quantization error for $\P$
of order $r$ is given by \cite{GL:00} 
\begin{eqnarray}
e_{n,r}(\P)\coloneqq\left\{ \begin{array}{ll}
\inf_{\alpha\in\mathcal{D}_{n}}\big(\int d(x,\alpha)^{r}\:\mathrm{d}\P(x)\big)^{1/r},\;\;\;\;\;\; r>0,\\
\inf_{\alpha\in\mathcal{D}_{n}}\exp\int\log d(x,\alpha)\:\mathrm{d}\P(x),\;\;\; r=0.
\end{array}\right.\label{quanerrordef}
\end{eqnarray}
According to \cite{GL:00}, $e_{n,r}(\P)$ equals the error with respect
to the $L_{r}$-metrics in the approximation of $\P$ with discrete
probability measures supported on at most $n$ points. See \cite{GL:00,KZQNumbers}
for various equivalent definitions for the quantization error. In
the following we will focus on the $L_{r}$-quantization problem with
$r>0$. For the quantization with respect to the geometric mean error,
we refer to \cite{GL:04} for rigorous foundations and \cite{Zhu:12b1,Zhu:Markov,key-11,key-17}
for more related results.

The upper and lower quantization dimension for $\P$ of order $r$,
as defined below, characterize the asymptotic quantization error in
a natural manner: 
\[
\overline{D}_{r}(\P)\coloneqq\limsup_{n\to\infty}\frac{\log n}{-\log e_{n,r}(\P)},\;\;\underline{D}_{r}(\P)\coloneqq\liminf_{n\to\infty}\frac{\log n}{-\log e_{n,r}(\P)}.
\]
If $\overline{D}_{r}(\P)=\underline{D}_{r}(\P)$, we call the common
value the quantization dimension of $\P$ of order $r$ and denote
it by $D_{r}(\P)$. To obtain more accurate information about the
asymptotic quantization error, we define the $s$-dimensional upper
and lower quantization coefficient (cf. \cite{GL:00,PK:01}): 
\[
\overline{Q}_{r}^{s}(\P)\coloneqq\limsup_{n\to\infty}n^{1/s}e_{n,r}(\P),\;\;\underline{Q}_{r}^{s}(\P)\coloneqq\liminf_{n\to\infty}n^{1/s}e_{n,r}(\P),\;\; s>0.
\]
By \cite{GL:00,PK:01}, the upper (lower) quantization dimension is
exactly the critical point at which the upper (lower) quantization
coefficient jumps from zero to infinity.

The following theorem by Zador is a classical result on quantization
of absolutely continuous measures. It was first proposed by Zador \cite{Za:63}
and then generalized by Bucklew and Wise \cite{BW:82}; we refer to \cite[Theorem 6.2]{GL:00}
for a rigorous proof. 
\begin{thm}
[\cite{GL:00}] Let $\P$ be absolutely continuous Borel probability
measure on $\mathbb{R}^{q}$ with density $h$ with respect to the
$q$-dimensional Lebesgue measure $\lambda^{q}$. Assume that for
some $\delta>0$, we have $\int|x|^{r+\delta}\:\mathrm{d}\P\left(x\right)<\infty$.
Then for all $r>0$ we have 
\[
\underline{Q}_{r}^{q}(\P)=\overline{Q}_{r}^{q}(\P)=C(r,q)\bigg(\int h^{\frac{q}{q+r}}(x)\:\mathrm{d}\lambda^{q}\left(x\right)\bigg)^{\frac{q+r}{q}},
\]
where $C(r,q)$ is a constant independent of $\P$. 
\end{thm}
While engineers are mainly dealing with absolutely continuous distributions,
the quantization problem is significant for all Borel probability
measures satisfying the moment condition $\int|x|^{r}\mathrm{d}\P\left(x\right)<\infty$.
For later use we define the subset of Borel probabilities $\mathcal{M}_{r}\coloneqq\left\{ \P\colon\P\left(\mathbb{R}\right)=1,\int|x|^{r}\mathrm{d}\P\left(x\right)<\infty\right\} $
and let $\mathcal{M}_{\infty}$ denote the set of Borel probability
measures with compact support. This condition ensures that the set of $n$-optimal
sets of order $r$ denoted by $C_{n,r}(\P)$ is non-empty. Also note
that $\mathcal{M}_{\infty}\subset\mathcal{M}_{r}$ for all $r>0.$
The most prominent aspects in quantization of probability measures
are the following:

\subsubsection*{Find the exact value of the upper/lower quantization dimension
for $\P$ of order $r$:}

In the case where the quantization dimension does not exist, it is
usually difficult to obtain the exact value of the upper or lower
one (cf. \cite{Roych:14}). Up to now, in such a situation, the
upper and lower quantization dimension could only be explicitly determined  for very special cases.

\subsubsection*{Determine the ${s}$-dimensional upper and lower quantization coefficient: }

We are mainly concerned about the finiteness and positivity of these
quantities. This question is analogous to the question of whether a fractal
is an $s$-set. Typically, this question is much harder to answer
than finding the quantization dimension. So far, the quantization
coefficient has been studied for absolutely continuous probability
measures (\cite{GL:00}) and several classes of singular measures,
including self-similar and self-conformal \cite{LM:02,Roych:11b,Zhu:11b,key-19} measures,
Markov-type measures \cite{kz:14Markov,Zhu:Markov,Roych:11b} and self-affine
measures on Bedford-McMullen carpets \cite{kz:13McMullen,Zhu:11}.

\subsubsection*{Properties of the point density measure ${\P_{r}}$:}

Fix a sequence of $n$-optimal sets $\alpha_{n}\in C_{n,r}(\P)$ of
order $r$, $n\in\mathbb{N}$, and consider the weak limit of the
empirical measures, whenever it exists, 
\[
\P_{r}\coloneqq\lim_{n\to\infty}\frac{1}{n}\sum_{a\in\alpha_{n}}\delta_{a}.
\]
The point density measure characterizes the frequency at which optimal
points fall into a given open set. Up to now, the point density measure
is determined only for absolutely continuous measures \cite[Theorem 7.5]{GL:00}
and certain self-similar measures \cite[Theorem 5.5]{GL:05}.

\subsubsection*{Local properties and Voronoi partitions: }

Fix a finite subset $\alpha$ of $\mathbb{R}^{q}$. A Voronoi partition
with respect to $\alpha$ refers to a partition $\left(P_{a}(\alpha)\right)_{a\in\alpha}$
of $\mathbb{R}^{q}$ such that 
\[
P_{a}(\alpha)\subset\left\{ x\in\mathbb{R}^{q}:d(x,\alpha)=d(x,a)\right\} ,\; a\in\alpha.
\]
It is natural to ask, if there exists constants $0<C_{1}\leq C_{2}<\infty$
such that for all $\alpha_{n}\in C_{n,r}(\P)$ and $n\in\mathbb{N}$
we have 
\begin{eqnarray*}
\frac{C_{1}e_{n,r}^{r}}{n} & \leq & \min_{a\in\alpha_{n}}\int_{P_{a}\left(\alpha_{n}\right)}\,\mathrm{d}(x,\alpha_{n})^{r}\:\mathrm{d}\P(x)\\
 & \leq & \max_{a\in\alpha_{n}}\int_{P_{a}\left(\alpha_{n}\right)}\,\mathrm{d}(x,\alpha_{n})^{r}\:\mathrm{d}\P(x)\leq\frac{C_{2}e_{n,r}^{r}}{n}.
\end{eqnarray*}
This question is essentially a weaker version of Gersho's conjecture
\cite{Ger:79}. Graf, Luschgy and Pag\`es proved in \cite{GL:12} that
this is in fact true for a large class of absolutely continuous probability
measures. An affirmative answer is also given for self-similar measures
under the assumption of the strong separation condition (SSC) for
the corresponding iterated function system \cite{Zhu:11b,Zhu:12b}.

In the final analysis, the study of the quantization problem addresses
the  optimal sets. Where do the points of an optimal set lie?
Unfortunately, it is almost impossible to determine the optimal sets
for a general probability measure. It is therefore desirable to seek
for an ``approximately explicit'' description of such sets. In other
words, even though we do not know exactly where the points of an optimal
set lie, we want to know how many points are lying in a given open
set. This would in return enable us to obtain precise estimates for
the quantization error.

\subsubsection*{Connection to fractal geometry: }

To this end, some typical techniques in fractal geometry are often
very helpful. In fact, the quantization problem is closely connected
with important notions in fractal geometry. One may compare the upper
(lower) quantization dimension for measures to the packing (Hausdorff)
dimension for sets; accordingly, the upper (lower) quantization coefficient
may be compared to the packing (Hausdorff) measure for sets. Although
they are substantially different, they do have some close connections,
as all these quantities can be defined in terms of coverings, partitions
and packings. In fact, we have 
\begin{enumerate}[labelindent=2pt,leftmargin=*,label=(\arabic*),widest=IV,align=left,font=\upshape\mdseries]
\item $\dim_{H}^{*}\mu\leq\underline{D}_{r}(\mu)\leq\underline{\dim}_{B}^{*}\mu$
and $\dim_{P}^{*}\mu\leq\overline{D_{r}}(\mu)\leq\overline{\dim}_{B}^{*}\mu$,
for $r=2$ these inequalities were presented in \cite{PK:01}, and
for measures with compact support and all $r\in(0,\infty]$ they were
independently proved in \cite{GL:00}. 
\item In \cite{KZStability} we have studied the \emph{stability} of the
upper and lower quantization dimension in some detail. In \cite{KZStability},
for $\ensuremath{r\in\left[1,\infty\right]}$, we proved the following
statments:

\begin{enumerate}[labelindent=2pt,leftmargin=*,label=(\roman*),widest=IV,align=left,font=\upshape\mdseries]
\item $\overline{D}_{r}\left(\mu\right)\coloneqq\inf\left\{ {\displaystyle \sup_{1\leq i\leq n}}\overline{D}_{r}\left(\mu_{i}\right):\mu=\sum_{i=1}^{n}s_{i}\mu_{i},\mu_{i}\in\mathcal{M}_{r},s_{i}>0,n\in\mathbb{N}\right\} $
for $\mu\in\mathcal{M}_{r}$. 
\item $\dim_{P}^{*}\left(\mu\right)\coloneqq\inf\left\{ {\displaystyle \sup_{i\in\mathbb{N}}}\overline{D}_{r}\left(\mu_{i}\right):\mu=\sum_{i\in\mathbb{N}}s_{i}\mu_{i},\mu_{i}\in\mathcal{M}_{\infty},s_{i}>0\right\} $
for all $\mu\in\mathcal{M}_{\infty}$. 
\item There exists $\mu\in\mathcal{M}_{\infty}$ such that $\overline{D}_{r}\left(\mu\right)\not=\dim_{P}^{*}\left(\mu\right)$. 
\item There exists $\mu\in\mathcal{M}_{\infty}$ such that $\overline{D}_{r}(\mu)>\underline{D}_{r}(\mu)$. 
\item There exists $\mu\in\mathcal{M}_{\infty}$ such that $\underline{D}_{r}\left(\mu\right)$
and its finitely stabilized counterpart $\inf\left\{ {\displaystyle \sup_{1\leq i\leq n}}\underline{D}_{r}\left(\mu_{i}\right):\mu=\sum_{i=1}^{n}s_{i}\mu_{i},\,\mu_{i}\in\mathcal{M}_{\infty},s_{i}>0,n\in\mathbb{N}\right\} $
\newline do not coincide. 
\end{enumerate}
\item for certain measures arising from dynamical systems, the quantization
dimension can be expressed within the thermodynamic formalism in terms
of appropriate\emph{ temperature functions }(see \cite{LM:02,kz:13McMullen,Roych:10,Roych:11a}). 
\item The upper and lower quantization dimension of order zero are closely
connected with the \emph{upper and lower local dimension}. As it is
shown in \cite{Zhu:12b}, if $\nu$-almost everywhere the upper and
lower local dimension are both equal to $s$, then $D_{0}(\nu)$ exists
and equals $s$. 
\end{enumerate}
We end this section with Graf and Luschgy's results on self-similar
measures. These results and the methods involved in their proofs have
a significant influence on subsequent work on the quantization for
non-self-similar measures.

Let $(S_{i})_{i=1}^{N}$ be a family of contractive similitudes on
$\mathbb{R}^{q}$ with contraction ratios $(s_{i})_{i=1}^{N}$. According
to \cite{Hut:81}, there exists a unique non-empty compact subset
$E$ of $\mathbb{R}^{q}$ such that $E=\bigcup_{i=1}^{N}S_{i}(E)$.
The set $E$ is called the self-similar set associated with $(S_{i})_{i=1}^{N}$.
Also, there exists a unique Borel probability measure on $\mathbb{R}^{q}$,
such that $\P=\sum_{i=1}^{N}p_{i}\P\circ S_{i}^{-1}$, called the
self-similar measure associated with $(S_{i})_{i=1}^{N}$ and the
probability vector $(p_{i})_{i=1}^{N}$. We say that $(S_{i})_{i=1}^{N}$
satisfies the\emph{ strong separation condition} (SSC) if the sets
$S_{i}(E),i=1,\cdots,N$, are pairwise disjoint. We say that it satisfies
the \emph{open set condition} (OSC) if there exists a non-empty open
set $U$ such that $S_{i}(U)\cap S_{j}(U)=\emptyset$ for all $i\neq j$
and $S_{i}(U)\subset U$ for all $i=1,\cdots,N$. For $r\in[0,\infty)$,
let $k_{r}$ be the positive real number given by 
\begin{eqnarray}
k_{0}\coloneqq\frac{\sum_{i=1}^{N}p_{i}\log p_{i}}{\sum_{i=1}^{N}p_{i}\log s_{i}},\;\;\sum_{i=1}^{N}(p_{i}c_{i}^{r})^{\frac{k_{r}}{k_{r}+r}}=1.\label{gl}
\end{eqnarray}

\begin{thm} [\cite{GL:01,GL:04}] Assume that $(S_{i})_{i=1}^{N}$
satisfies the open set condition. Then for all $r\in[0,\infty)$,
we have 
\[
0<\underline{Q}_{r}^{k_{r}}(\P)\leq\overline{Q}_{r}^{k_{r}}(\P)<\infty.
\]
In particular, we have $D_{r}(\P)=k_{r}$. \end{thm} 

This is the first complete result on the quantization for (typically)
singular measures. In its proof, H\"older's inequality with an exponent
less than one plays a crucial role, from which the exponent $k_{r}/\left(k_{r}+r\right)$
comes out in a natural manner.

\section{The three-step procedure}\label{3step}

Following the ideas of Graf-Luschgy we propose a three-step procedure
for the estimation of the quantization errors by means of partitions,
coverings and packings. This procedure is applicable to a large class
of fractal measures, including Moran measures, self-affine measures
and Markov-type measures, provided that some suitable separation condition
is satisfied; it even allows us to obtain useful information on the
quantization for general Borel probability measures on $\mathbb{R}^{q}$
with compact support.

\subsubsection*{{Step~1 (Partitioning)} }

For each $n$, we partition the (compact) support of $\mu$ into $\varphi_{n}$
small parts $\left(F_{nk}\right)_{k=1}^{\varphi_{n}}$, such that
$\mu(F_{nk})|F_{nk}|^{r}$ are uniformly comparable, namely, for some
constant $C>1$ independent of $k,j\in\left\{ 1,\ldots,\varphi_{n}\right\} $
and $n\in\mathbb{N}$, we have 
\[
C^{-1}\mu(F_{nk})|F_{nk}|^{r}\leq\mu(F_{nj})|F_{nj}|^{r}\leq C\mu(F_{nk})|F_{nk}|^{r},
\]
where $|A|$ denotes the diameter of a set $A\subset\mathbb{R}^{d}$.
This idea was first used by Graf and Luschgy to treat the quantization
problem for self-similar measures, we refer to \cite{GL:00} for a
construction of this type. The underlying idea is to seek for some
uniformity while $\mu$ generally is not uniform.

\subsubsection*{Step~2 (Covering) }

With a suitable separation condition, we may also assume that for
some $\delta>0$, we have that 
\[
d(F_{nk},F_{nj})\geq\delta\max\{|F_{nk}|,|F_{nj}|\},\; k\neq j,\; n\geq1.
\]
In this step, uniformity and separation allow us to verify that any
$\varphi_{n}$-optimal set distributes its points equally among suitable
neighborhoods of $F_{nk},$ $1\leq k\leq\varphi_{n}$, in other words,
each $F_{nk}$ ``owns'' a bounded number of points of the $\varphi_{n}$-optimal
set. More precisely, we prove that there exists some constant $L_{1}$,
which is independent of $n$, such that for every $\alpha\in C_{\varphi_{n},r}(\mu)$,
we have 
\[
\max_{1\leq k\leq\varphi_{n}}{\rm card}\left(\alpha\cap(F_{nk})_{4^{-1}\delta|F_{nk}|}\right)\leq L_{1},
\]
where $A_{s}$ denotes the $s$\emph{-parallel set} of $A$. This
can often be done inductively by means of contradiction.

\subsubsection*{Step~3 (Packing)}

In the last step we have to find a constant $L_{2}$ and subsets $\beta_{nk}$ of $F_{nk}$
with cardinality at most  $L_{2}$ such that for all   $\alpha\in C_{\varphi_{n},r}(\mu)$ and $x\in F_{nk}$ we have
\[
d(x,\alpha)\geq d((x,(\alpha\cap(F_{nk})_{4^{-1}\delta|F_{nk}|})\cup\beta_{nk}).
\]
This reduces the global situation to a local one and enables us to
restrict our attention to an arbitrary small set $F_{nk}$. We have
\[
e_{\phi_{n},r}^{r}(\mu)\geq\sum_{k=1}^{\phi_{n}}\int_{F_{nk}}d(x,(\alpha\cap(F_{nk})_{4^{-1}\delta|F_{nk}|})\cup\beta_{nk})^{r}\,\mathrm{d}\mu(x).
\]
Note that ${\rm card}\left((\alpha\cap(F_{nk})_{4^{-1}\delta|F_{nk}|})\cup\beta_{nk}\right)\leq L_{1}+L_{2}$.
For measures with explicit mass distributions, we often have 
\[
\int_{F_{nk}}d(x,\gamma\cup\beta_{k})^{r}\,\mathrm{d}\mu(x)\geq D\mu(F_{nj})|F_{nj}|^{r}
\]
for any subset $\gamma$ of $\mathbb{R}^{q}$ with cardinality not
greater than $L_{1}+L_{2}$ and an appropriate constant $D$. Thus, we get a lower estimate for the
quantization error: 
\[
e_{\phi_{n},r}^{r}(\mu)\geq D\sum_{k=1}^{\phi_{n}}\mu(F_{nk})|F_{nk}|^{r}.
\]
On the other hand, by choosing some arbitrary points $b_{k}\in F_{nk}$,
$k\in\left\{ 1,\ldots,\phi_{n}\right\} $, one can easily see 
\[
e_{\phi_{n},r}^{r}(\mu)\leq\sum_{k=1}^{\phi_{n}}\int_{F_{nk}}d(x,b_{k})^{r}d\mu(x)\leq\sum_{k=1}^{\phi_{n}}\mu(F_{nk})|F_{nk}|^{r}.
\]

After these three steps, for sufficiently ``nice'' measures, we may
additionally assume that $\phi_{n}\leq\phi_{n+1}\leq C\phi_{n}$ for some constant
$C>1$ (cf. \cite{key-22,key-23,key-20}). To determine the dimension it
is then enough to estimate the growth rate of $\phi_{n}$. Here, 
ideas from Thermodynamic Formalism --~such as critical exponents
or zeros of some pressure function~-- often come into play: E.g., for $r>0$
 we often have 
\[
\frac{D_{r}(\P)}{D_{r}(\P)+r}=\inf\bigg\{ t\in\mathbb{R}:\sum_{n\in\mathbb{N}}\sum_{k=1}^{\phi_{n}}\left(\mu(F_{nk})|F_{nk}|^{r}\right)^{t}<\infty\bigg\}
\]
allowing us to find explicit formulae for the quantization dimension
for a given problem (see \cite{kz:13McMullen} for an instance of
this). Typically, for a non-self-similar
measure such as a self-affine measures on Bedford-McMullen carpets,
this requires a detailed analysis of the asymptotic quantization errors. 
 In order to formulate a rigorous proof,
we usually need to make more effort according to the particular properties
of the measures under consideration. As general measures do not enjoy
strict self-similarity, it seems unrealistic to expect to establish
simple quantities for the quantization errors as Graf and Luschgy
did for self-similar measures \cite[Lemma 14.10]{GL:00}. However,
the above-mentioned three-step procedure often provides us with estimates
of the quantization errors which is usually a promising starting point.

Moreover, in order to examine the finiteness or positivity of the
upper and lower $s$-dimensional quantization coefficient of order
$r$, it suffices to check that (cf. \cite{Zhu:12a}) 
\[
0<\liminf_{n\to\infty}\sum_{k=1}^{\phi_{n}}(\mu(F_{nj})|F_{nj}|^{r})^{\frac{s}{s+r}}\leq\limsup_{n\to\infty}\sum_{k=1}^{\phi_{n}}(\mu(F_{nj})|F_{nj}|^{r})^{\frac{s}{s+r}}<\infty.
\]
An effective way to do this is to construct some auxiliary probability
measures. Such a measure should closely reflect the information carried
by $(\mu(F_{nj})|F_{nj}|^{r})^{\frac{s}{s+r}}$. For a self-similar
measure, as Graf-Luschgy's work shows, an auxiliary probability measure
is the self-similar measure associated with $(S_{i})_{i=1}^{N}$ and
the probability vector $((p_{i}c_{i}^{r})^{\frac{k_{r}}{k_{r}+r}})_{i=1}^{N}$.
It is interesting to note that this measure coincides with the point
density measure provided that the $k_r$-dimensional quantization coefficient exists.
For a self-similar measure, as Graf and Luschgy showed, we can use
the above auxiliary probability measure and obtain the finiteness
or positivity of the upper and lower $k_{r}$-dimensional quantization
coefficient, which also implies that the quantization dimension exists
and equals $k_{r}$.  In the non-self-similar situation, due to the complexity of the
topological support, it is often not easy to construct a suitable
auxiliary probability measure to estimate the quantization coefficients. 

\section{Recent work on the quantization for fractal measures}

\subsection{Self-affine measures on Bedford-McMullen carpets}

Fix two positive integers $m,n$ with $2\leq m\leq n$ and fix a set
\[
G\subset\big\{0,1,\ldots,n-1\big\}\times\big\{0,1,\ldots,m-1\big\}
\]
with $N\coloneqq\mbox{card}\left(G\right)\geq2$. We define a family
of affine mappings on $\mathbb{R}^{2}$ by 
\begin{equation}
f_{ij}:(x,y)\mapsto\big(n^{-1}x+n^{-1}i,m^{-1}y+m^{-1}j\big),\;\;(i,j)\in G.\label{fi's}
\end{equation}
By \cite{Hut:81}, there exists a unique non-empty compact set $E$
satisfying $E=\bigcup_{(i,j)\in G}f_{ij}(E)$, which is called the
Bedford-McMullen carpet determined by $(f_{ij})_{(i,j)\in G}$. Given
a probability vector $(p_{ij})_{(i,j)\in G}$ with $p_{ij}>0$, for
all $(i,j)\in G$, the self-affine measure associated with $(p_{ij})_{(i,j)\in G}$
and $(f_{ij})_{(i,j)\in G}$ refers to the unique Borel probability
measure $\mu$ on $\mathbb{R}^{2}$ satisfying 
\begin{equation}
\mu=\sum_{(i,j)\in G}p_{ij}\mu\circ f_{ij}^{-1}.\label{affine}
\end{equation}
Sets and measures of this form have been intensively studied in the
past decades, see e.g. \cite{Bed:84,Mcmullen:84,LG:92,Peres:94b,King:95,Fal:10,GLi:10}
for many interesting results. We write 
\begin{eqnarray*}
G_{x} & \coloneqq & \left\{ i:(i,j)\in G\;\mbox{for some }j\right\} ,\; G_{y}\coloneqq\left\{ j:(i,j)\in G\;\mbox{for some }i\right\} ,\\
G_{x,j} & \coloneqq & \left\{ i:(i,j)\in G\right\} ,\;\; q_{j}\coloneqq\sum_{i\in G_{x,j}}p_{ij}.
\end{eqnarray*}
We carry out the three-step procedure and obtain an estimate for the
quantization errors. This allows us to conjecture
that the quantization dimension exists and equals $s_{r}$, where
\begin{eqnarray}
\bigg(\sum_{(i,j)\in G}(p_{ij}m^{-r})^{\frac{s_{r}}{s_{r}+r}}\bigg)^{\theta}\bigg(\sum_{j\in G_{y}}(q_{j}m^{-r})^{\frac{s_{r}}{s_{r}+r}}\bigg)^{1-\theta}=1,\;\theta\coloneqq\frac{\log m}{\log n}.\label{maineq1}
\end{eqnarray}
However, it seems rather difficult to find a suitable auxiliary measure
for a proof of this conjecture. A cornerstone is the crucial observation
that the number $s_{r}$ coincide with a Poincare-like exponent \cite{kz:13McMullen}.
Using the property of sup-additive sequences, we are able to prove
that $D_{r}(\mu)$ exists and also coincides with $\kappa_{r}$. Finally,
we consider the self-affine measure associated with $((p_{ij}m^{-r})^{\frac{s_{r}}{s_{r}+r}}/C_{r})_{(i,j)\in G}$
as an auxiliary measure, where $C_{r}\coloneqq\sum_{(i,j)\in G}(p_{ij}m^{-r})^{\frac{s_{r}}{s_{r}+r}}$.
This measure and the above-mentioned estimate enable us to obtain
sufficient conditions for the upper and lower quantization coefficient
to be both positive and finite. We have 

\begin{thm} [\cite{kz:13McMullen}] Let $\mu$ be as defined in
(\ref{affine}). Then for each $r\in(0,\infty)$ we have that $D_{r}(\mu)$
exists and equals $s_{r}$, Moreover, $0<\underline{Q}_{r}^{s_{r}}(\mu)\leq\overline{Q}_{r}^{s_{r}}(\mu)<\infty$
if one of the following conditions is fulfilled: 
\begin{enumerate}[labelindent=\parindent,leftmargin=*,label=(\Alph*),widest=IV,align=left,font=\upshape\mdseries]
\item \label{ConditionA} $\sum_{i\in G_{x,j}}(p_{ij}q_{j}^{-1})^{\frac{s_{r}}{s_{r}+r}}$
are identical for all $j\in G_{y}$, 
\item \label{ConditionB}$q_{j}$ are identical for all $j\in G_{y}$. 
\end{enumerate}
\end{thm} \textbf{Open problem}: Is it true that $0<\underline{Q}_{r}^{s_{r}}(\mu)\leq\overline{Q}_{r}^{s_{r}}(\mu)<\infty$
if and only if condition \ref{ConditionA} or \ref{ConditionB} holds?

\subsection{Quantization for Markov-type measures}

\subsubsection{Mauldin-Williams fractals}

Let $J_{i}$, non-empty compact subsets of $\mathbb{R}^{d}$ with
$J_{i}={\rm cl}({\rm int}(J_{i}))$, $1\leq i\leq N$, where ${\rm cl}(A)$ and ${\rm int}(A)$
denote the closure and interior in $\mathbb{R}^{d}$ of a set $A\subset\mathbb{R}^{d}$.
 For the integer $N\geq 2$ let $P=(p_{ij})_{1\leq i,j\leq  N}$
be a row-stochastic matrix, i.e., $p_{ij}\geq0,1\leq i,j\leq N$,
and $\sum_{j=1}^{N}p_{ij}=1,1\leq i\leq N$. Let $\theta$ denote
the empty word and set
\begin{eqnarray*}
\Omega_{0} & \coloneqq & \{\theta\},\;\Omega_{1}\coloneqq\{1,\ldots N\},\\
\Omega_{k} & \coloneqq & \{\sigma\in\Omega_{1}^{k}:p_{\sigma_{1}\sigma_{2}}\cdots p_{\sigma_{k-1}\sigma_{k}}>0\},\; k\geq2,\\
\Omega^{*} & \coloneqq & \bigcup_{k\geq0}\Omega_{k},\;\Omega_{\infty}\coloneqq\{\sigma\in\Omega_{1}^{\mathbb{N}}:p_{\sigma_{h}\sigma_{h+1}}>0\;\;{\rm for\; all}\;\; h\geq1\}.
\end{eqnarray*}
We call $J_{i},1\leq i\leq N$, cylinder sets of order one. For each
$1\leq i\leq N$, let $J_{ij},\;(i,j)\in\Omega_{2}$, be non-overlapping
subsets of $J_{i}$ such that $J_{ij}$ is geometrically similar to
$J_{j}$ and ${\rm diam}(J_{ij})/{\rm diam}(J_{j})=c_{ij}$. We call
these sets cylinder sets of order two. Assume that cylinder sets of
order $k$ are determined, namely, for each $\sigma\in\Omega_{k}$,
we have a cylinder set $J_{\sigma}$. Let $J_{\sigma\ast i_{k+1}},\sigma\ast i_{k+1}\in\Omega_{k+1}$,
be non-overlapping subsets of $J_{\sigma}$ such that $J_{\sigma\ast i_{k+1}}$
is geometrically similar to $J_{i_{k+1}}$. Inductively, cylinder
sets of order $k$ are determined for all $k\geq1$. The {\em (ratio specified)
Mauldin-Williams fractal} is given by 
\[
E\coloneqq\bigcap_{k\geq1}\bigcup_{\sigma\in\Omega_{k}}J_{\sigma}.
\]

\subsubsection{Markov-type measures}

Let $(\chi_{i})_{i=1}^{N}$ be an arbitrary probability vector with
$\min_{1\leq i\leq N}\chi_{i}>0$. By Kolmogorov consistency theorem,
there exists a unique probability measure $\widetilde{\mu}$ on $\Omega_{\infty}$
such that $\widetilde{\mu}([\sigma])\coloneqq\chi_{\sigma_{1}}p_{\sigma_{1}\sigma_{2}}\cdots p_{\sigma_{k-1}\sigma_{k}}$
for every $k\geq1$ and $\sigma=(\sigma_{1},\ldots,\sigma_{k})\in\Omega_{k}$,
where $[\sigma]\coloneqq\{\omega\in\Omega_{\infty}:\omega|_{|\sigma|}=\sigma\}$.
Let $\pi$ denote the projection from $\Omega_{\infty}$ to $E$ given
by $\pi\left(\sigma\right)\coloneqq x$, where 
\[
\{x\}\coloneqq\bigcap_{k\geq1}J_{\sigma|_{k}},\;\;{\rm for}\;\;\sigma\in\Omega_{\infty}.
\]
Let us assume the following:
\begin{enumerate}[labelindent=\parindent,leftmargin=*,label=(A\arabic*),widest=IV,align=left,font=\upshape\mdseries]
\item  ${\rm card}(\{j:p_{ij}>0\})\geq2$ for all $1\leq i\leq N$.
\item There exists a constant $t\in\left(0,1\right)$ such that for every
$\sigma\in\Omega^{*}$ and distinct $i_{1},i_{2}\in\Omega_{1}$ with
$\sigma\ast i_{l}\in\Omega_{\left|\sigma\right|+1}$, $l=1,2$, 
\[
d(J_{\sigma\ast i_{1}},J_{\sigma\ast i_{2}})\geq t\max\{|J_{\sigma\ast i_{1}}|,|J_{\sigma\ast i_{2}}|\}.
\]

\end{enumerate}
Under this assumption, $\pi$ is a bijection. We consider the image
measure of $\widetilde{\mu}$ under the projection $\pi$ given by
$\mu\coloneqq\widetilde{\mu}\circ\pi^{-1}$. We call $\mu$ a Markov-type
measure which satisfies 
\begin{eqnarray}
\mu(J_{\sigma})=\chi_{\sigma_{1}}p_{\sigma_{1}\sigma_{2}}\cdots p_{\sigma_{k-1}\sigma_{k}}\;\;{\rm for}\;\;\sigma=(\sigma_{1}\ldots\sigma_{k})\in\Omega_{k}.\label{markovmeasure}
\end{eqnarray}

For $1\leq i,j\leq N$, we define $a_{ij}(s)\coloneqq(p_{ij}c_{ij}^{r})^{s}$.
Then we get an $N\times N$ matrix $A(s)=(a_{ij}(s))_{N\times N}$.
Let $\psi(s)$ denote the spectral radius of $A(s)$. By \cite[Theorem 2]{MW:88},
$\psi(s)$ is continuous and strictly decreasing. Note that, by the
assumption (A1), the Perron-Frobenius theorem and intermediate-value
theorem, there exists a unique number $\xi\in(0,1)$ such that $\psi(\xi)=1$.
Thus, for every $r>0$, there exists a unique positive number $s_{r}$
such that $\psi(\frac{s_{r}}{s_{r}+r})=1$.

We consider the directed graph $G$ associated with the transition
matrix $(p_{ij})_{N\times N}$. Namely, $G$ has vertices $1,2,\ldots,N$.
There is an edge from $i$ to $j$ if and only if $p_{ij}>0$. In
the following, we will simply denote by $G=\{1,\ldots,N\}$ both the
directed graph and its vertex sets. We also write 
\[
b_{ij}(s)\coloneqq(p_{ij}c_{ij}^{r})^{\frac{s}{s+r}},\;\; A_{G,s}\coloneqq(b_{ij}(s))_{N\times N},\;\;\Psi_{G}(s)\coloneqq\psi\bigg(\frac{s}{s+r}\bigg).
\]
Let ${\rm SC}(G)$ denote the set of all strongly connected components
of $G$. For $H_{1},H_{2}\in{\rm SC}(G)$, we write $H_{1}\prec H_{2}$,
if there is a path initiating at some $i_{1}\in H_{1}$ and terminating
at some $i_{k}\in H_{2}$. If we have neither $H_{1}\prec H_{2}$
nor $H_{2}\prec H_{1}$, then we say $H_{1},H_{2}$ are incomparable.
For every $H\in{\rm SC}(G)$, we denote by $A_{H,s}$ the sub-matrix
$(b_{ij}(s))_{i,j\in H}$ of $A_{G}(s)$. Let $\Psi_{H}(s)$ be the
spectral radius of $A_{H,s}$ and $s_{r}(H)$ be the unique positive
number satisfying $\Psi_{H}(s_{r}(H))=1$.

Again, we apply the three-step procedure in section 2 and obtain upper
and lower estimates for the quantization error. Using these estimates
and auxiliary measures of Mauldin-Williams type, we are able to prove
that, when the transition matrix is irreducible, the upper and lower
quantization coefficient are both positive and finite. This fact also
leads to the positivity of the lower quantization coefficient in the general
case. Then, based on a detailed analysis of the corresponding directed
graph (not strongly connected) and some techniques in matrix theory,
we are able to prove the formula for the quantization dimension. Finally,
by using auxiliary measures of Mauldin-Williams type once more,
we establish a necessary and sufficient condition for the upper quantization
coefficient to be finite as stated next. 

\begin{thm} [\cite{kz:14Markov}] Assume that (A1) and (A2) are
satisfied. Let $\mu$ be the Markov-type measure as defined in (\ref{markovmeasure})
and $s_{r}$ the unique positive number satisfying $\Psi_{G}(s_{r})=1$.
Then, $D_{r}(\mu)=s_{r}$ and $\underline{Q}_{r}^{s_{r}}(\mu)>0$.
Furthermore, $\overline{Q}_{r}^{s_{r}}(\mu)<\infty$ if and only if
$\mathcal{M}\coloneqq\{H\in{\rm SC}(G):s_{r}(H)=s_{r}\}$ consists
of incomparable elements, otherwise, we have $\underline{Q}_{r}^{s_{r}}(\mu)=\infty$.
\end{thm}

\subsection{Quantization for Moran measures}

\subsubsection{Moran sets}

Let $J$ be a non-empty compact subset of $\mathbb{R}^{d}$ with $J={\rm cl}({\rm int}(J))$.
Let $|A|$ denote the diameter of a set $A\subset\mathbb{R}^{d}$.
Let $(n_{k})_{k=1}^{\infty}$ be a sequence of integers with $\min_{k\geq1}n_{k}\geq2$
and $\theta$ denote the empty word. Set 
\[
\Omega_{0}\coloneqq\{\theta\},\;\Omega_{k}\coloneqq\prod_{j=1}^{k}\{1,2,\cdots,n_{j}\},\;\Omega^{*}\coloneqq\bigcup_{k=0}^{\infty}\Omega_{k}.
\]
For $\sigma=\sigma_{1}\cdots\sigma_{k}\in\Omega_{k}$ and $j\in\{1,\cdots,n_{k+1}\}$,
we write $\sigma\ast j=\sigma_{1}\cdots\sigma_{k}j$.

Set $J_{\theta}\coloneqq J$ and let $J_{\sigma}$ for $\sigma\in\Omega_{1}$
be non-overlapping subsets of $J_{\theta}$ such that each of them
is geometrically similar to $J_{\theta}$. Assume that $J_{\sigma}$
is determined for every $\sigma\in\Omega_{k}$. Let $J_{\sigma\ast j},1\leq j\leq n_{k+1}$
be non-overlapping subsets of $J_{\sigma}$ which are geometrically
similar to $J_{\sigma}$. Inductively, all sets $J_{\sigma},\sigma\in\Omega^{*}$
are determined in this way. The Moran set is then defined by 
\begin{equation}
E\coloneqq\bigcap_{k=1}^{\infty}\bigcup_{\sigma\in\Omega_{k}}J_{\sigma}.\label{moransets}
\end{equation}
We call $J_{\sigma},\sigma\in\Omega_{k}$, cylinders of order $k$.
It is well known that the Moran sets $E$ are generally not self-similar
(cf. \cite{CM:93,Wen:01}). For $k\geq0$ and $\sigma\in\Omega_{k}$,
we set 
\[
|\sigma|\coloneqq k,\; c_{\sigma,j}\coloneqq\frac{|J_{\sigma\ast j}|}{|J_{\sigma}|},\;1\leq j\leq n_{k+1}.
\]
We assume that there exist some constants $c,\beta\in(0,1)$ such
that 
\begin{enumerate}[labelindent=\parindent,leftmargin=*,label=(B\arabic*),widest=IV,align=left,font=\upshape\mdseries]
\item ${\displaystyle \inf_{\sigma\in\Omega^{*}}\min_{1\leq j\leq n_{|\sigma|+1}}c_{\sigma,j}=c>0}$,
\label{enu:B1}
\item ${\rm dist}(J_{\sigma\ast i},J_{\sigma\ast j})\geq\beta\max\{|J_{\sigma\ast i}|,|J_{\sigma\ast j}|\}$
for $1\leq i\neq j\leq n_{|\sigma|+1}$ and $\sigma\in\Omega^{*}$.
\label{enu:B2}
\end{enumerate}

\subsubsection{Moran measures}

For each $k\geq1$, let $(p_{kj})_{j=1}^{n_{k}}$ be a probability
vector. By the Kolmogorov consistency theorem, there exists a probability
measure $\nu_{\omega}$ on $\Omega_{\infty}:=\prod_{k=1}^{\infty}\{1,2,\cdots,n_{k}\}$
such that 
\begin{eqnarray*}
\nu([\sigma_{1},\cdots,\sigma_{k}])=p_{1\sigma_{1}}\cdots p_{k\sigma_{k}},\;\sigma_{1}\cdots\sigma_{k}\in\Omega_{k},
\end{eqnarray*}
where $[\sigma_{1},\cdots,\sigma_{k}]=\{\tau\in\Omega_{\infty}:\tau_{j}=\sigma_{j},1\leq j\leq k\}$.
Let $\Pi:\Omega_{\infty}\to E$ be defined by $\Pi(\sigma)=\bigcap_{k\geq1}J_{\sigma|_{k}}$
with $\sigma|_{k}=\sigma_{1}\cdots\sigma_{k}$. Then, with the assumption
\ref{enu:B2}, $\Pi$ is a continuous bijection. We define $\mu:=\nu\circ\Pi^{-1}$.
Then, we have 
\[
\mu(J)=1,\;\mu(J_{\sigma}):=p_{1\sigma_{1}}\cdots p_{k\sigma_{k}},\;\;\sigma=\sigma_{1}\cdots\sigma_{k}\in\Omega_{k},\; k\geq1.
\]
We call the measure $\mu$ the Moran measure on $E$. It is known
that the quantization dimension for $\mu$ of order $r$ does not
necessarily exist. Let $d_{k,r},\overline{d}_{r},\underline{d}_{r}$
be given by 
\begin{eqnarray*}
\sum_{\sigma\in\Omega_{k}}(p_{\sigma}c_{\sigma}^{r})^{\frac{d_{k,r}}{d_{k,r}+r}}=1,\;\;\overline{d}_{r}:=\limsup_{k\to\infty}d_{k,r},\; d_{k,r}\underline{d}_{r}:=\liminf_{k\to\infty}d_{k,r}.
\end{eqnarray*}
\textbf{Open problem} Is it true that $\overline{D}_{r}(\mu)=\overline{d}_{r},\underline{D}_{r}(\mu)=\underline{d}_{r}$?

%\subsection{Multiscale Moran measures}

\subsubsection{Multiscale Moran sets}

A multiscale Moran set is Moran set with some additional structure
encoded in the infinite sequence $\omega=(\omega_{l})_{l=1}^{\infty}\in\Upsilon\coloneqq\{1,\dots,m\}^{\mathbb{N}}$
for some $m\geq2$. For this fix some positive integers $N_{i}\geq2,1\leq i\leq m$
and for every $1\leq i\leq m$, let $(g_{ij})_{j=1}^{N_{i}}$ be the
contraction vector with $g_{ij}\in(0,1)$ and $(p_{ij})_{j=1}^{N_{i}}$
a probability vector with $p_{ij}>0$ for all $1\leq j\leq N_{i}$.
Now using the notation in the definition of Moran sets, we set 
\begin{eqnarray}
n_{l+1}\coloneqq N_{\omega_{l+1}},\;\;\left(c_{\sigma,j}\right)_{j=1}^{N_{\omega_{l+1}}}\coloneqq(g_{\omega_{l+1}j})_{j=1}^{N_{\omega_{l+1}}},\;\sigma\in\Omega_{l},\; l\geq0.\label{fi2}
\end{eqnarray}
If, for some $l\geq0$, we have $\omega_{l+1}=i$, then for every
$\sigma\in\Omega_{l}$, we have a continuum of choices of $\{J_{\sigma\ast j}\}_{j=1}^{N_{i}}$
fulfilling \ref{enu:B1}, \ref{enu:B2} and (\ref{fi2}), because
we only fix the contraction ratios of the similitudes. Hence, to every
$\omega\in\Upsilon$, there corresponds a class $\mathcal{M}_{\omega}$
of Moran sets according to (\ref{moransets}). We call these Moran
sets \emph{multiscale Moran sets}.

For each $\omega\in\Upsilon$, we write 
\[
N_{k,i}(\omega)\coloneqq{\rm card}\{1\leq l\leq k:\omega_{l}=i\},\;\;1\leq i\leq m.
\]
Fix a probability vector $\chi=(\chi_{i})_{i=1}^{m}$ with $\chi_{i}>0$
for all $1\leq i\leq m$ and define 
\begin{eqnarray*}
G(\chi) & \coloneqq & \{\omega\in\Upsilon:\lim_{k\to\infty}k^{-1}N_{k,i}(\omega)=\chi_{i},\;1\leq i\leq m\},\\
G_{0}(\chi) & \coloneqq & \{\omega\in\Upsilon:\;\limsup_{k\to\infty}\big|N_{k,i}(\omega)-k\chi_{i}\big|<\infty,1\leq i\leq m\}.
\end{eqnarray*}

\subsubsection{Multiscale Moran measures}

Fix an $\omega\in G(\chi)$. According to Kolmogorov consistency theorem,
there exists a probability measure $\nu_{\omega}$ on the product
space $\Omega_{\infty}\coloneqq\prod_{k=1}^{\infty}\{1,2,\cdots,N_{\omega_{k}}\}$
such that 
\begin{eqnarray*}
\nu_{\omega}([\sigma_{1},\cdots,\sigma_{k}])=p_{\omega_{1}\sigma_{1}}\cdots p_{\omega_{k}\sigma_{k}},\;\sigma_{1}\cdots\sigma_{k}\in\Omega_{k},
\end{eqnarray*}
where $[\sigma_{1},\cdots,\sigma_{k}]=\{\tau\in\Omega_{\infty}:\tau_{j}=\sigma_{j},1\leq j\leq k\}$.
We define $\mu_{\omega}\coloneqq\nu_{\omega}\circ\Pi^{-1}$. Then,
we have 
\[
\mu(J)=1,\;\mu_{\omega}(J_{\sigma})\coloneqq p_{\omega_{1}\sigma_{1}}\cdots p_{\omega_{k}\sigma_{k}},\;\;\sigma=\sigma_{1}\cdots\sigma_{k}\in\Omega_{k},\; k\geq1.
\]
We call the measure $\mu_{\omega}$ the infinite product measure on
$E(\omega)$ associated with $\omega$ and $(p_{ij})_{j=1}^{N_{i}},1\leq i\leq m$.

For every $\omega\in G(\chi)$ and $k\in\mathbb{N}$, let $s_{k,r}(\omega),s_{r}$
and $\underline{H}_{r}(\omega),\overline{H}_{r}(\omega)$ be defined
by 
\begin{eqnarray}
 &  & \prod_{i=1}^{m}\bigg(\sum_{j=1}^{N_{i}}(p_{ij}g_{ij}^{r})^{\frac{s_{k,r}(\omega)}{s_{k,r}(\omega)+r}}\bigg)^{N_{k,i}(\omega)}=1,\;\;\prod_{i=1}^{m}\bigg(\sum_{j=1}^{N_{i}}(p_{ij}g_{ij}^{r})^{\frac{s_{r}}{s_{r}+r}}\bigg)^{\chi_{i}}=1,\label{kz1}\\
 &  & \underline{H}_{r}(\omega)\coloneqq\liminf_{k\to\infty}k|s_{k,r}(\omega)-s_{r}|,\;\;\overline{H}_{r}(\omega)\coloneqq\limsup_{k\to\infty}k|s_{k,r}(\omega)-s_{r}|.\nonumber 
\end{eqnarray}
Compared with Mauldin-Williams fractals, the disadvantage is that
we have more patterns in the construction of multiscale Moran sets.
However, the pattern we use at the $(k+1)$-th step is independent
of words of length $k$, which is an advantage. After we carry out
the three-step procedure in Section \ref{3step}, we conveniently obtain the exact
value of the quantization dimension by considering some measure-like
auxiliary functions. This also enables us to transfer the question
of the upper and lower quantization coefficient to the convergence
order of $(s_{k,r}(\omega))_{k=1}^{\infty}$. For the latter, we need
a detailed analysis of some auxiliary functions related to (\ref{kz1}).
One may see \cite{Zhu:12a} for more details. Our main result is summarized
in the following theorem. 

\begin{thm} [\cite{Zhu:12a}] For every $\omega\in G(\chi)$, we
have 
\begin{enumerate}[labelindent=\parindent,leftmargin=*,label=(\roman*),widest=IV,align=left,font=\upshape\mdseries]
\item $D_{r}(\mu_{\omega})$ exists and equals $s_{r}$, it is independent
of $\omega\in G(\chi)$, 
\item if $s_{k,r}(\omega)\geq s_{r}$ for all large $k$, then $\underline{Q}_{r}^{s_{r}}(\mu_{\omega})>0$.
If in addition $\underline{H}_{r}(\omega)=\infty$, then we have $\overline{Q}_{r}^{s_{r}}(\mu_{\omega})=\infty$,
\item if $s_{k,r}(\omega)\leq s_{r}$ for all large $k$, then $\overline{Q}_{r}^{s_{r}}(\mu_{\omega})<\infty$;
if, in addition, $\underline{H}_{r}(\omega)=\infty$, then we have
$\underline{Q}_{r}^{s_{r}}(\mu_{\omega})=0$, 
\item \label{enu:Hr<infty}if $\overline{H}_{r}(\omega)<\infty$, then $\underline{Q}_{r}^{s_{r}}(\mu_{\omega})$
and $\overline{Q}_{r}^{s_{r}}(\mu_{\omega})$ are both positive and
finite,
\item if $\omega\in G_{0}(\chi)$, then the assertion in \emph{\ref{enu:Hr<infty}}
holds. 
\end{enumerate}
\end{thm} 

\textbf{Open problem:} What can we say about necessary conditions
for $\underline{Q}_{r}^{s_{r}}(\mu_{\omega})$ and $\overline{Q}_{r}^{s_{r}}(\mu_{\omega})$
to be both positive and finite?

\end{document}